\documentclass{article}
\usepackage[margin=1in]{geometry}
\usepackage{comment}
\usepackage{float}
\usepackage{graphicx}
\usepackage{authblk}
\usepackage[normalem]{ulem}
\usepackage{amssymb}
\usepackage{amsmath}

\title{Data-driven discovery of governing equations for coarse-grained heterogeneous network dynamics}

\author[]{Katherine Owens}
\author[]{J. Nathan Kutz}
\affil[]{\em Department of Applied Mathematics, University of Washington, Seattle WA}

\begin{document}

\maketitle

\hrule
\vspace{.1in}
\noindent
{\bf Abstract:}
We leverage data-driven model discovery methods to determine the governing equations for the emergent behavior of heterogeneous networked dynamical systems. Specifically, we consider networks of coupled nonlinear oscillators whose collective behaviour approaches a limit cycle. Stable limit-cycles are of interest in many biological applications as they model self-sustained oscillations (e.g. heart beats, chemical oscillations, neurons firing, circadian rhythm).  For systems that display relaxation oscillations, our method automatically detects boundary (time) layer structures in the dynamics, fitting inner and outer solutions and matching them in a data-driven manner. We demonstrate the method on well-studied systems: the Rayleigh Oscillator and the Van der Pol Oscillator. We then apply the mathematical framework to discovering low-dimensional dynamics in networks of semi-synchronized Kuramoto, Rayleigh, Rossler, and Fitzhugh-Nagumo oscillators, as well as heterogeneous combinations thereof. We also provide a numerical exploration of the dimension of collective network dynamics as a function of several network parameters, showing that the discovery of coarse-grained variables and dynamics can be accomplished with the proposed architecture.   
\vspace{.1in}
\hrule

\section{Introduction}
\begin{figure}[t]
\includegraphics[width=\textwidth]{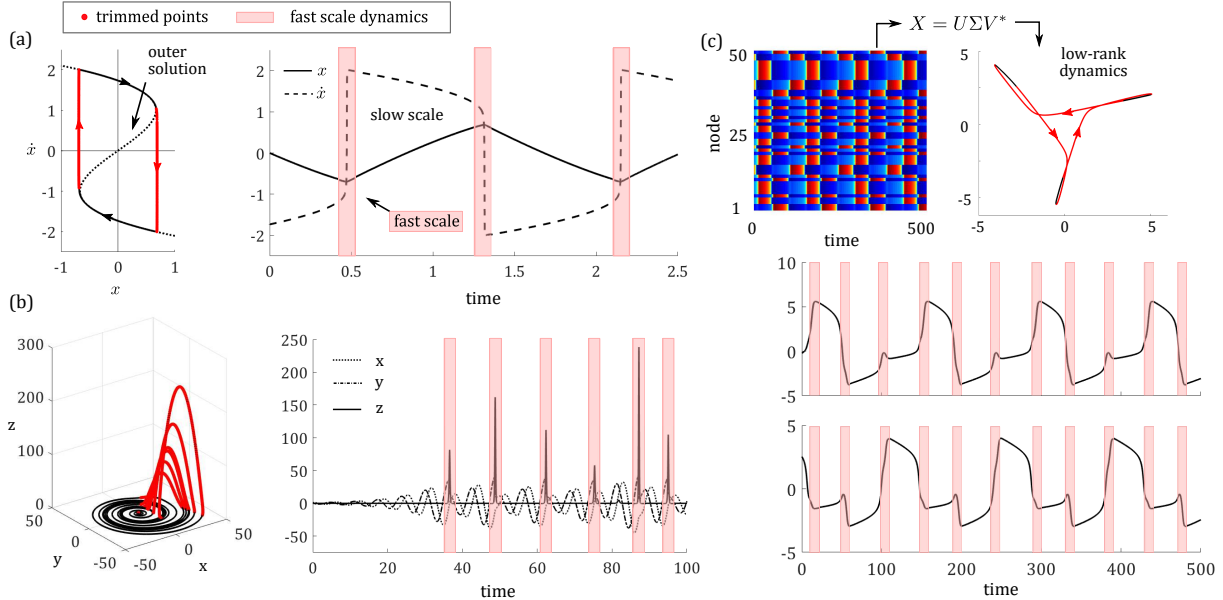}
\caption{Sparse Identification of Nonlinear Dynamics (SINDy) with data-trimming disambiguates fast-scale versus slow-scale dynamics. (a) For a sample trajectory on the limit cycle of a Rayleigh oscillator, the trimmed points are those that deviate from the asymptotic outer solution. (b) For a sample trajectory on the Rossler attractor, trimming isolates rapid excursions in the z plane. (c) The dynamics of a network of coupled Fitzhugh-Nagumo oscillators can be projected onto the first 2 principal components. In this low dimensional representation, SINDy with trimming isolates intervals of spiking and rapid relaxation. }
\label{fig: trimming}
\end{figure}

\begin{figure}[t]
\includegraphics[width=\textwidth]{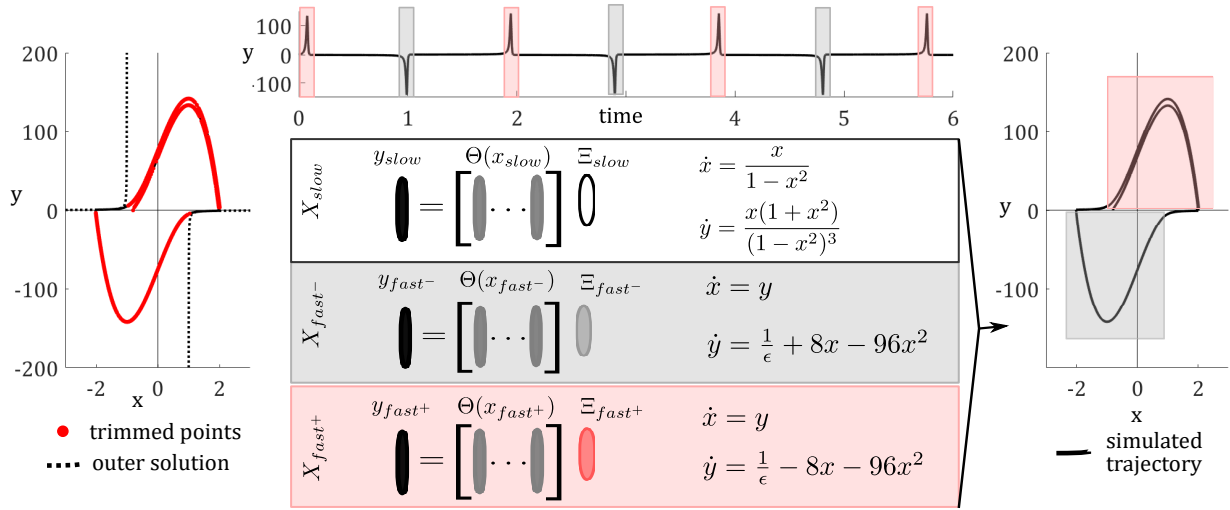}
\caption{SINDy with trimming for the Van der Pol Oscillator. As with the Rayleigh example, data is partitioned by performing SINDy with trimming. The trimmed points (red) are those that deviate from the outer solution (dotted line). Then, SINDy is applied to the points from each time scale successively. Models are fit to the fast scale dynamics in the lower half-plane and the upper half-plane separately, but we can see that under the change of variable $\hat{x} = -x$, the resulting equations for $\dot{y}$ in the gray and red boxes are the same. Combining these equations in the appropriate regions, the hybrid model can then be used to simulate the Van der Pol system.}
\label{fig: VanderPol}
\end{figure}

Data-driven modeling paradigms are an increasingly important tool for the characterization of behavior in complex, high-dimensional systems. For instance, {\em  Reduced order models} (ROMs)~\cite{benner2015survey} aim to exploit the ubiquitous observation from experiment and/or computation that meaningful input/output signals in high-dimensional systems are encoded in low-dimensional patterns of dynamic activity. However, many complex systems of interest exhibit behavior across multiple temporal or spatial scales, which poses unique challenges for modeling and predicting their behavior. It is often the case that while we are primarily interested in macroscale phenomena, the microscale dynamics must also be modeled and understood, as they play a central role in driving the macroscale behavior. This can make dealing with multiscale systems particularly difficult unless the time scales are disambiguated in a principled way. There is a significant body of research focused on modeling multiscale systems to produce coarse-grained descriptions: notably the {\em equation-free} (EF) method\cite{kevrekidis_equation-free_2003} for linking scales and the {\em heterogeneous multiscale modeling} (HMM) framework  \cite{weinan_heterognous_2003}. Additional work has focused on testing for the presence of multiscale dynamics so that analyzing and simulating multiscale systems is more computationally efficient \cite{froyland_computational_2014,froyland_trajectory-free_2016}. Many of the same issues that make modeling multiscale systems difficult can also present challenges for model discovery and system identification. These challenges motivate the development of specialized methods for performing model discovery on problems with multiple time scales, taking into account the unique properties of multiscale systems. Specifically, we integrate dimensionality-reduction, sparse regression, and robust statistics to discover governing evolution equations for coarse-grained heterogeneous network dynamics as well as to explicitly disambiguate and model multiscale temporal phenomena (See Fig. \ref{fig: VanderPol}).

The discovery and pairing of coordinates and dynamics is the hallmark feature of scientific modeling.  For instance, the second century Ptolemaic description of celestial mechanics featured an earth-centric coordinate systems with the retrograde planetary dynamics expressed as linear superposition of circular motion of different frequencies, i.e. the doctrine of the perfect circle.  By Newton's time in the seventeenth century, the heliocentric coordinate system was established allowing for an ${\bf F}=m{\bf a}$ description of planetary motion with the inverse square law of gravitation.  Thus it is not surprising that methods for jointly finding coordinates and dynamics have emerged as a central theme in modern data-driven science.  Indeed, there is significant diversity in mathematical approaches for pairing a coordinate representation with an accompanying dynamical model that characterizes the evolution. {\em Reduced order models} (ROMs) accomplish this pairing by constructing coordinates from the most dominant correlated activity, or {\em proper orthogonal decomposition} (POD)~\cite{holmes2012turbulence,benner2015survey}, and then Galerkin projecting on the governing equations~\cite{benner2015survey}. More broadly, modal decomposition techniques, such as POD and {\em dynamic mode decomposition} (DMD)~\cite{kutz2016dynamic}, approximate {\em linear} subspaces using dominant correlations in spatio-temporal data~\cite{Taira2017aiaa}.  Linear subspaces, however, are highly restrictive and ill-suited to handle parametric dependencies. Attempts to circumvent these shortcomings include using multiple linear subspaces covering different temporal or spatial domains \cite{Dihlmann11,bright2013compressive,Taddei2015,Peherstorfer2015}, multi-resolution DMD~\cite{kutz2016multiresolution},
diffusion map embeddings~~\cite{coifman2005pnas,coifman2006acha,coifman2008mmas,nadler2006acha}, or more recently, using deep learning to compute underlying nonlinear subspaces which are advantageous for dynamics, both linear and nonlinear~\cite{brunton2016discovering,lusch2018deep,champion2019data,champion2019discovery,Pan2020}. These techniques represent data-driven architectures for extracting order-parameter descriptions of the underlying spatio-temporal dynamics observed~\cite{cross1993pattern}.

Such flexibility and diversity in algorithms for the joint discovery of coordinates and dynamics has allowed for many new representations of physical systems, including in the networked dynamical systems of interest here. Specifically, we are interested in modeling the emergent coarse-grained heterogenous behavior that arises in high-dimensional, networked dynamical systems whose underlying dynamics are strongly nonlinear. In particular, we consider networks of coupled nonlinear oscillators whose collective behavior approaches relaxation-oscillation limit cycle dynamics. Recently, oscillatory phenomena have been coarse-grained into learned and/or collective coordinates (intrinsic coordinates) for modeling their evolution~\cite{kemeth2020learning,gottwald2015model,cartwright2019collective,shlizerman2012neural,morrison2020nonlinear,nishikawa2003heterogeneity}. 
Gottwald originally introduced collective coordinates to describe coupled phase oscillators of the Kuramoto model~\cite{gottwald2015model}, and then extended the mathematical framework to reduce infinite-dimensional stochastic partial differential equations (SPDEs) with symmetry to a set of finite-dimensional stochastic differential equations which describe the shape of the solution and the dynamics along the symmetry group~\cite{cartwright2019collective}. Shlizerman et al~\cite{shlizerman2012neural} coarse-grained using coordinates associated with the dominant correlated activity, or singular value decomposition, and projected into a framework of neural activity measures. Kemeth et al~\cite{kemeth2020learning} instead learn collective dynamics on a slow manifold, after initial transients have died out, which can be approximated through a learned model based on local {\em spatial} partial derivatives in the emergent coordinates. These works are closely related to the present work, as their aims are well-aligned with the goals of our method. However, there are key differences. Specifically, we consider coarse-graining on a system for which regions of fast activity are critical for characterizing the dynamics, i.e. the dynamics produce temporal boundary layers.  In particular, we wish to learn the dominant balance dynamics that occur in both fast and slow temporal regimes on network level dynamics.

In this work, we leverage dimensionality-reduction, sparse regression, and robust statistics to discover coarse-grained models of heterogeneous networked dynamical systems. We consider networks consisting of several types of oscillators and explore how changing the network structure influences the collective dynamics.
Nishikawa et al~\cite{nishikawa2003heterogeneity} and Motter et al~\cite{motter2005network} consider synchronization of networks, including scale free and small world, as a function of the homogeneous or heterogeneous distribution of connectivity, showing that homogeneous networks are more easily synchronized.  Here, we focus on aspects of the emergent dynamics that arise from the heterogeneity of the oscillators themselves.  Specifically, in many regimes, the emergent dynamics exhibit two key features we aim to model: (i) nontrivial coarse-grained dynamics (relaxation-oscillations), and (ii) temporal boundary layers (rapid transient dynamics). Specifically, the proposed mathematical architecture discovers parsimonious governing evolution equations for coarse-grained network dynamics and, if desired, explicitly disambiguates emergent temporal boundary layer phenomena. The developed algorithm starts simply by projecting onto coordinates associated with the dominant correlated activity, i.e. by taking the singular value decomposition. We then identify rapid transition regions associated with transients in the time dynamics of the dominant modes by applying the {\em sparse identification of nonlinear dynamics} (SINDy))~\cite{brunton2016discovering} algorithm with data-trimming~\cite{champion2020unified}, which is a standard analysis tool from robust statistics. The governing equations produced by this stage of the method can be used to describe the full limit cycle. However, if a more detailed model is desired, the results of trimming are used to segment the coarse-grained trajectory, separating regions of rapid transition from the slow-temporal field, which is equivalent to the slow center-manifold considered previously~\cite{kemeth2020learning,gottwald2015model}. We then utilize the SINDy framework again to discover parsimonious governing equations for each region of the coarse-grained trajectory, constructing a hybrid model for the overall dynamics. Here we apply the full method to a network of coupled Fitzhugh-Nagumo oscillators, which produce rapid (spiking) behavior. The strongly nonlinear nature of the underlying agents leads to the existence of temporal boundary layers in the collective dynamics. The overall architecture allows for a data-driven approach to characterizing dominant-balance physics~\cite{cross1993pattern,callaham2021learning} and singularly perturbed problems~\cite{bender2013advanced,kevorkian2013perturbation,kutz2020advanced}, which is a classic technique for understanding, for  instance, boundary layer formation in fluid dynamics. 

\begin{figure}[t]
\includegraphics[width=\textwidth]{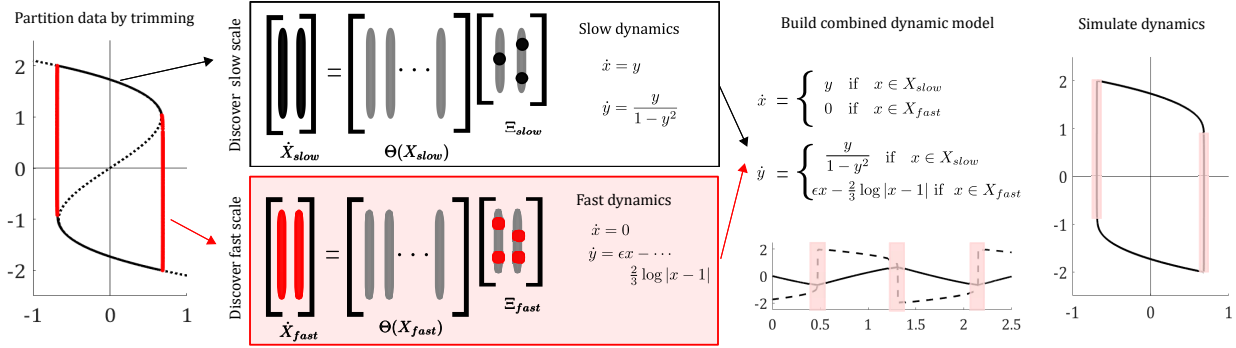}
\caption{Overview of SINDy with trimming for discovering fast and slow dynamics illustrated on the Rayleigh Oscillator. First, data is partitioned by performing SINDy with trimming. Next, SINDy is applied to the points from each time scale separately. Subsequently, a combined dynamic model is built which applies the fast or slow model appropriately at each step depending on the state of the system. Simulated dynamics closely match the input trajectory.}
\label{fig: Rayleigh}
\end{figure}

\section{Coarse-graining and SINDy for multi-scale systems}

Coarse-graining is of broad scientific interest across a diversity of disciplines. Coarse-graining seeks to determine a set of variables that successfully summarize the detailed state of a network of interacting agents.  Analytical and computational approaches have long been sought for achieving this end, with the former leveraging the dynamical systems theory of center manifold reductions, normal forms, bifurcation theory and perturbation analysis.  Indeed, the seminal review of Cross and Hohenberg~\cite{cross1993pattern} detail the many mathematical architectures available for deriving order-parameter descriptions, or coarse-grained models, that characterize the dynamics.  More recent computational approaches include statistical and machine learning algorithms~\cite{bishop2006pattern,murphy2012machine} for network clustering, principal component analysis, and, for instance, algorithms for detecting density (k-core decomposition).  The diversity of methods developed allow for a new coordinate representation for the collective dynamics. This representation can then be leveraged by model discovery methods such as SINDy to determine the dynamics on the coarse-grained coordinates.

SINDy recovers parsimonious representations of the dynamics from measurement data by sparse regression to a library of candidate models~\cite{brunton2016discovering}.  Thus the goal of SINDy is to promote a parsimonious representation of the complex system by traditional and interpretable governing equations. 
Consider the nonlinear dynamical system, which will ultimately be our coarse-grained representation,
\begin{equation}
    \frac{d{\bf x}}{dt} = f({\bf x},t)
  \label{eq:plato}
\end{equation}
which is measured at time points ${\bf t}=[t_1, t_2, \cdots, t_m]$. From measurements, we construct the matrix ${\bf X}({\bf t})$ = $[{{\bf x}_1({\bf t}) \,\, {\bf x}_2({\bf t})\,\,  \dots \,\, {\bf x}_n({\bf t})}]$ $\in \mathbb{R}^{m\times n}$. The method introduced in ~\cite{brunton2016discovering} seeks to identify ${\bf f}$ via sequential threshold least-squares, which is a proxy for the sparsifying zero-norm. The set of $n$ state measurements are used to populate a library of candidate nonlinear terms ${\bf \Theta (X) = [1^\top \; X^\top \; (X \otimes X)^\top \; \cdots \; \text{sin}(X)^\top ]}$, where ${\bf x \otimes y}$ defines the vector of all product combinations of the state components. Each candidate term should be unique, as a suitable library is crucial in the SINDy algorithm. A common strategy is start with polynomials and increase the complexity of the library with other terms, such as trigonometric functions. Thus, the system in (\ref{eq:plato}) is approximated by:
\begin{equation}
    {\bf \dot{X} = \Theta(X)\Xi } .
  \label{eq:SINDy}
\end{equation}
The time derivatives ${\bf \dot{X}}({\bf t}) = [ {\bf \dot{x}}_1({\bf t}) \,\, {\bf \dot{x}}_2({\bf t}) \,\, \dots \,\, {\bf \dot{x}}_n({\bf t})]$, if not measured directly, can be found via numerical differentiation and should be appropriately de-noised, if necessary (low pass Butterworth, total variation regularization, etc.)~\cite{van2020numerical}. The coefficients ${\bf \Xi}$ are the {\em sparse} weightings of the corresponding candidate library terms. Therefore, our regression relies on sparse regularization to enforce a parsimonious $\bf \Xi$ corresponding to the fewest nonlinear terms in our library that describe our dynamics well: 
\begin{equation}
    {\bf \Xi} = \text{arg}\,\min\limits_{\hat{ \bf \Xi}\, } \Vert{\bf \Theta(X)\hat{\bf \Xi} - \dot{X}}\Vert_\text{2} + \lambda\Vert{\bf \hat{\bf \Xi}}\Vert_\text{0}
  \label{eq:sparse}
\end{equation}
Regressing to the zero-norm is often achieved by relaxing the the one-norm.  However, modern optimization frameworks are allowing for computationally tractable proxies for the zero-norm that are superior to the one-norm relaxation~\cite{champion2020unified,zheng2018unified}.
SINDy is modular and adaptible, allowing for modifications which include the discovery of spatio-temporal systems~\cite{rudy2017data,schaeffer2017learning}, multiscale dynamics~\cite{champion2019data}, parametric dependencies~\cite{rudy2019data}, hybrid dynamical systems~\cite{mangan2019model}, systems subject to control~\cite{kaiser2018sparse,brunton2016sparse}, and implicit dynamics~\cite{mangan2016inferring,mangan2017model,kaheman2020sindy}.  Further, it can be used in the low-data limit~\cite{quade2018sparse,hirsh2021sparsifying,fasel2021ensemble} and with denoising architectures~\cite{rudy2019smoothing,kaheman2020automatic}.  The SINDy algorithm, which has an accompanying python package~\cite{kaptanoglu2021pysindy}, is the workhorse algorithm for our coarse-graining method.

To disambiguate fast scale dynamics, our coarse-graining method exploits a key extension of SINDy: data trimming. Modifying Eq. (\ref{eq:sparse}) to include trimming yields:
\begin{equation}
\begin{split}
 {\bf \Xi}, {\bf v} ~= ~& \text{arg}\,\min\limits_{\hat{ \bf \Xi}, {\bf v} } \sum_{i=1}^m \frac{1}{2} v_i \Vert{(\bf \Theta(X)\hat{\bf \Xi} - \dot{X})_i}\Vert_\text{2} + \lambda\Vert{\bf \hat{\bf \Xi}}\Vert_\text{0} \\
& \text{s.t.}~~~~ 0 \leq v_i \leq 1, ~~~~~ {\bf1}^T\bf{v} = h,
 \end{split}
    \label{eqn: trimming}
\end{equation}
where $h < m$ is an estimate of the number of the $m$ input data points that are considered ``inliers." See the work of Champion et al. for a thorough explanation \cite{champion2020unified}. After the SINDy with trimming algorithm is applied, the vector ${\bf v}$ has lower entries where data points are harder to fit to a parsimonious model of system dynamics. In systems that exhibit multi-scale temporal dynamics, these troublesome points often correspond to spikes in the derivative. Thus, trimming effectively identifies regions of rapid change. 

As an example, in Fig. \ref{fig: trimming}a we illustrate how data-trimming partitions the limit cycle of a Rayleigh oscillator: the points which deviate from the slow-scale outer solution are trimmed, coinciding with the two boundary layers in the limit cycle. Trimming applied to Rossler dynamics on a strange attractor identifies rapid excursions in the z-component of the system (Fig. \ref{fig: trimming}b). When applied to a low-rank representation of a network of spiking Fitzhugh-Nagumo (FHN) oscillators, trimming identifies regions of rapid spiking and relaxation (Fig. \ref{fig: trimming}c). In each of these examples, the trimmed fraction, or $1-h/m$, is a hyperparameter of the method that must be tuned according to the proportion of the input measurements that appear to belong to regions of fast-scale dynamics. If the fraction is set too high, then points adjacent to the fast scale dynamics will also be trimmed. Conversely, if the trimmed fraction is set too low then only a subset of the fast-scale dynamics will be identified. Another factor to consider in applying this method is that in order to effectively repurpose an outlier-detection technique from robust statistics to partition our data, the data must be clean. If the input trajectory is noisy, then data trimming will pull out the noisy data points and their neighbors, whose derivative estimates are corrupted, and, possibly, also the desired regions of rapid scale dynamics.

\section{Hybrid Models for Relaxation Oscillators}
 
The steps of our method for data-driven identification of temporal boundary layer phenomena in relaxation oscillations are outlined on the canonical Rayleigh oscillator in Fig. \ref{fig: Rayleigh}. First, SINDy with trimming is applied to the input trajectory to partition the data. Any reasonable library can be used at this point, as we only need the output from trimming, i.e. the vector ${\bf v}$ from (\ref{eqn: trimming}), not the model coefficients. The entries of ${\bf v}$ are compared to a threshold value, and points with entries below the threshold are assigned to the fast-scale group, $X_{fast}$, while points with an entry above the threshold are assigned to the slow-scale group, $X_{slow}$. Next, the locations of the points belonging to $X_{fast}$ within the phase space are used to bound the region in which the fast scale dynamics will apply. The bounding box is formed around each stretch of consecutive trimmed points, and then overlapping boxes are joined resulting in one region for each segment of the limit cycle displaying fast-scale dynamics. Points outside of these bounding boxes are assigned to $X_{slow}$ during simulation. Then, SINDy is applied separately to data in $X_{fast}$ and $X_{slow}$. For the slow-scale group, the candidate library of functions is built using $x \in X_{slow}$. For the fast scale group, the candidate library of functions is built using $x \in X_{fast}$. With these two dynamic models in hand, and a partition on the phase space determining when each model applies, we can simulate the system dynamics. At each time step, before calculating an update to the system position, we check whether the point belongs to $X_{fast}$ or $X_{slow}$. If $x \in X_{fast}$, then the fast scale dynamics are used to calculate an update. If $x \in X_{slow}$, then the slow scale dynamics are used to calculate an update. Note that the hybrid model is only valid for points on the limit cycle of the system. In this example, trajectories initiated elsewhere in the phase space either remain on the line $y = 0$, approach the line $y  = 1$ in the upper half plane, or approach the line $y = -1$ in the negative half plane. 

We also demonstrate this hybrid model approach on the limit cycle of a Van der Pol oscillator exhibiting relaxation-oscillations, as shown in Fig. \ref{fig: VanderPol}. Like the Rayleigh Oscillator example, we first partition the input trajectory using trimming. Note that the trimmed points, shown in red, are those points that deviate from the outer solution predicted by asymptotic approximation theory, shown by a dashed line. We then use the set of trimmed points to bound the regions of phase space governed by fast scale dynamics, and the remaining phase space is considered to be governed by slow scale dynamics. Next we apply SINDy to data points from the slow scale scale region, fitting the dynamics to a library that depends on $x$ and includes rational functions.  In this case, we fit the dynamics of the two regions of fast scale dynamics separately, calling those in the lower half plane $X_{fast^-}$ and those in the upper half plane $X_{fast^+}$. For both subsets of data, we fit the dynamics to a library of cubic polynomials in $x$ and $y$. We treat the two segments of fast dynamics separately because fitting both regions at once results in an average between the two desired models that does not approximate either branch fully. We note, however, that the resulting models for $X_{fast^-}$ and $X_{fast^+}$ reported in Fig. \ref{fig: VanderPol}, can easily be transformed to match each other by applying an absolute value sign to $x$. If such symmetries are clear from viewing the limit cycle of a system, then the library of candidate functions used in SINDy could be built using polynomials of $|x|$, and both regions of fast scale dynamics could be fit together. Again simulation of the system is accomplished by checking whether the current state belongs to  $X_{slow}, X_{fast^+},$ or $X_{fast^-}$, making an update with the appropriate dynamic rules, and then iterating these two steps until the desired time span is covered. 

\section{Coarse-graining heterogeneous spatio-temporal dynamics}
\begin{figure}[htp]
\includegraphics[width=\textwidth]{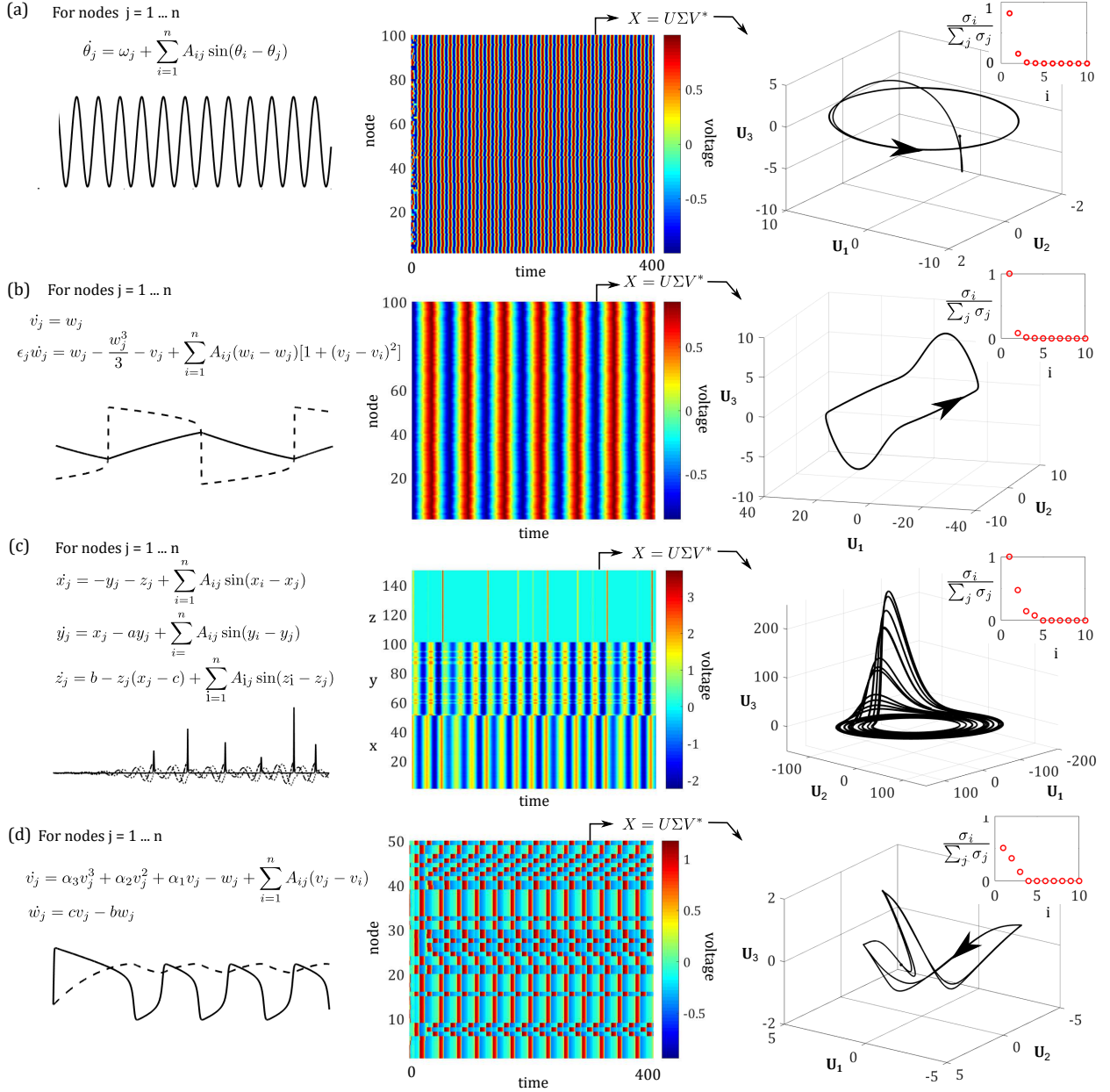}
\caption{In appropriate parameter regimes, the dynamics of networked oscillators can be reduced to a low dimensional representation by taking the singular value decomposition and projecting the state of the system at each time point onto the first few singular vectors. Here we illustrate the low-dimensional limit cycle that results for (a) coupled Kuramoto oscillators with different intrinsic frequencies, (b) coupled Rayleigh oscillators with different damping parameters, (c) linearly coupled Rossler oscillators, and (d) linearly coupled Fitzhugh-Nagumo (FHN) oscillators.}
\label{fig: coarsegraining}
\end{figure}
\begin{figure}[htp]
\includegraphics[width=\textwidth]{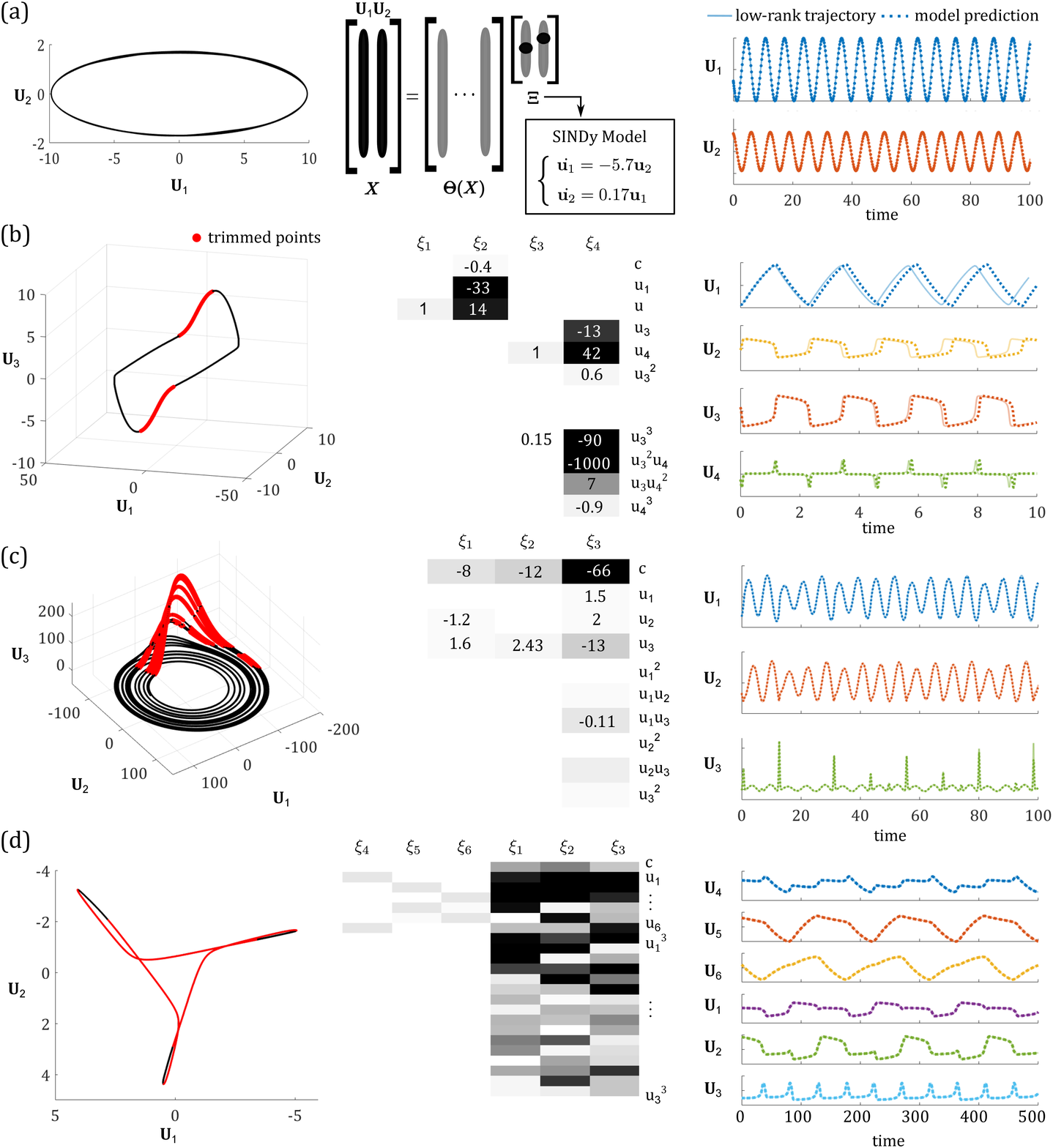}
\caption{SINDy applied to the low-dimensional representations of network dynamics produces governing equations which can be used to accurately simulate the low-rank dynamics over an appropriate time horizon. (a) For a network of Kuramoto oscillators, SINDy produces a very simple model that remains accurate for long time. (b) For a network of Rayleigh oscillators, SINDy produces a 3rd order model that comes out of phase after around one lap of the limit cycle. The discovered dynamics are attracted to a slightly smaller stable limit cycle than the actual system. (c) For a network of Rossler oscillators, SINDy produces a model that is linear in two dimensions and second order in the third. The predicted dyanamics sync up well for a while, but the prediction of spikes in $U_3$ starts to degrade, impacting the predictions of all three components over time. (d) For a network of FHN oscillators, trimming identifies regions of the limit cycle that correspond with semi-synchronized spiking. The resulting model is linear in the recovery-variable modes and 3rd order in the voltage variable modes. Although it is not sparse, the model predicts system dynamics accurately over a long time horizon.}
\label{fig: NetworksSINDy}
\end{figure}

For our coarse-graining method, a reduced-order basis is first derived through the singular value decomposition~\cite{kutz2013data,brunton2022data} (Fig. \ref{fig: coarsegraining}) and second, SINDy is deployed to discover the dynamics in this reduced order basis (Fig. \ref{fig: NetworksSINDy}). If the low-rank dynamics exhibit multiscale temporal behavior, then SINDy with trimming can be used instead to identify regions of fast-scale dynamics and fit a hybrid model to the data, facilitating analysis of the dominant balance physics in distinct dynamic intervals of the collective limit cycle. This method allows us to tackle coarse-graining for networks displaying heterogenous spatio-temporal dynamics that are not tractable with analytical approaches.  For very high-dimensional problems, randomized algorithms~\cite{erichson2016randomized} and compressed decompositions~\cite{erichson2019compressed} can be exploited in order to compute the low-rank subspaces more efficiently. 

The systems that we consider here are networks of coupled oscillators of the form 
\begin{equation}
    \frac{d{\bf x}}{dt} = f({\bf x},t),
\end{equation}
where ${{\bf x}(t)} = (x_1, x_2, .., x_n)$ is a vector of state variables, or nodes, and $f({\bf x} ) = (f_1({\bf x} ), f_2({\bf x }), ..., f_n({\bf x} ))$ is a vector of functions encoding the dynamics on each node. The dynamics depend on intrinsic properties of each node as well as coupling to neighboring nodes. Neighbors are defined by a symmetric, unweighted adjacency matrix, $A \in \mathbb{R}^{n \times n}$, where the entries $A_{ij} = A_{ji} = 1$ if node $i$ and node $j$ are connected and $A_{ij} = A_{ji} = 0$ otherwise. Here we consider Erdos-Renyi adjacency matrices, so $A$ may be sparse or highly connected as we tune the edge probability, $p$. The functions $f_i({\bf x})$ could take on any arbitrary form, but we limit our attention to a selection of well-known oscillators. 

Before applying our method to heterogeneous systems, we establish that it is producing reasonable results by considering a network consisting of Kuramoto oscillators with sinusoidal coupling. For the ubiquitous Kuramoto oscillator, the phase of node $j$, denoted $\theta_j$, is governed by 
\begin{equation}
\label{eqn: Kuramoto}
    \dot{\theta}_j = \omega_j + \frac{K}{n}\sum_{i=1}^n A_{ij}\sin(\theta_i - \theta_j).
\end{equation}
The key parameters besides the aforementioned adjacency matrix, are the connectivity strength  $K > 0$, and the intrinsic frequency $\omega_j $. For each node in the system, $\omega_j $ is drawn from a uniform distribution. Kuramoto oscillators have been widely studied and exhibit a well-known phase transition from asynchrony below a critical coupling threshold to  partial synchrony above this threshold. The system continues to synchronize more fully as the coupling parameter is further increased. Panel (a) of Figure \ref{fig: coarsegraining} illustrates the behavior of a single Kuramoto oscillator as well as the full system dynamics for a network of Kuramoto oscillators. Note that we are visualising the cosine of the phase. Taking the singular value decomposition and projecting onto the first few modes reduces the system dynamics to a stable elliptical limit cycle. Retaining two of these modes and applying the SINDy method to the 2D limit cycle produces a simple and accurate set of governing equations which can be used to simulate the system dynamics (Fig. \ref{fig: NetworksSINDy}a, final column). With this friendly, linear collective behavior no trimming is required to recover useful model equations.

Next, we consider a slightly more complicated system: Rayleigh oscillators with nonlinear coupling. In this case each node follows the second order governing equation
\begin{equation}
\label{eqn: Rayleigh}
    \epsilon_j \ddot{x}_j = \dot{x}_j - \frac{\dot{x}^3}{3} - x + \frac{K}{n}\sum_{i=1}^n A_{ij}(1+(x_j - x_i)^2)(\dot{x}_i-\dot{x}_j)
\end{equation}
where $\epsilon_j << 1$  is drawn from a uniform distribution and again $K>0$ is the connectivity parameter. The nonlinear coupling strategy, termed Haken-Kelso-Bunz (HKB) coupling, was introduced to study human bimanual experiments and applied to data capturing the synchronization of people rocking chairs by Alderisio et al. \cite{alderisio2016entrainment}. In numerical simulations, we observe partial synchrony in these systems when there is sufficient coupling (Fig. \ref{fig: coarsegraining}b). This collective behavior is captured by a nonlinear limit cycle in the low-dimensional representation. Applying SINDy to the 4D projection with trimming pulls out  regions of rapid relaxation along the direction of $U_4$ and produces a model that is linear in $U_1$ and $U_2$ and cubic in $U_3$ and $U_4$. Here each oscillator has an intrinsic $\epsilon_j$ on the order of $10^{-3}$. The multi-scale nature of the collective system dynamics is evident in the fact that the coefficients learned for modes $U_2$ and $U_4$ are much higher than those learned for $U_1$ and $U_3$. With a distribution of $\epsilon$ values centered closer to $1$, we would see a smaller gap between scales as $1/\epsilon \rightarrow 1$. Comparing the trajectory predicted by the model against the true low-rank trajectory of the system (final column of Fig. \ref{fig: NetworksSINDy}), we see that after one lap of the limit cycle, the prediction gets out of phase with the true trajectory. The SINDy model is attracted to a smaller stable limit cycle parallel to the true dynamics. Applying this method to systems with larger values of $\epsilon$, i.e. systems for which the boundary layers are less sharp, results in model dynamics that stay on track for longer time. 

We also demonstrate our method on a network of Rossler oscillators with sinusoidal coupling. In this case each node has 3 state variables governed by
\begin{equation}
\begin{split}
\label{eqn: Rossler}
\dot{x}_j & = -y_j - z_j + \frac{K}{n}\sum_{i=1}^n A_{ij}\sin(x_j - x_i) \\
\dot{y}_j & =  x_j + ay_j + \frac{K}{n}\sum_{i=1}^n A_{ij}\sin(y_j - y_i) \\
\dot{z}_j & = b + z_j(x_j-c) + \frac{K}{n}\sum_{i=1}^n A_{ij}\sin(z_j - z_i) 
\end{split}
\end{equation}
The parameters $(a,b,c)$ are uniform across all nodes in the network for a given run. In the example shown in Fig. \ref{fig: coarsegraining}c, $(a = 0.2, b = 0.2, c = 5.7)$.  The trajectory of each node is drawn to a similarly shaped strange attractor, however the attractors may be translated in phase space depending on the initial condition of the node. The nodes of the network largely synchronize phases and frequencies in their $x$ component, synchronize frequencies of their $y$ component, and partially synchronize the timing of excursions in the $z$ component. Projecting onto the first three singular vectors, the system lands on an attractor that resembles a Rossler attractor with a few differences. Specifically, the attractor is tilted slightly with respect to the axes, so there is a background oscillation rate in $U_3$. Due to this tilt, the $U_3$ ``excursions" also disrupt the background activity along $U_1$ and $U_2$ slightly. Applying SINDy with trimming to this low-dimensional trajectory results in a linear model for $U_1$ and $U_2$ and a dense cubic equation for the dynamics of $U_3$. Because the attractor for the collective dynamics is tilted in the SVD-space compared to an attractor for a single Rossler oscillator, these cross terms are necessary to capture the background frequency and the excursions. Observe that trimmed points are mostly spikes along the $U_3$  direction. In this case, the model predictions and the true trajectory match up very well for time steps that were used in the SVD projection initially. However, about time  $80$ in Fig. \ref{fig: NetworksSINDy}c you can start to see that the model prediction is deviating from the true trajectory--capturing less of the $U_3$ excursions. Continuing to predict beyond the training data, the model still produces dynamics that resemble a collective Rossler oscillator but they lose fidelity to the true system dynamics. 

The last type of oscillator that we considered were linearly coupled Fitzhugh-Nagumo oscillators. This simplified model of a neuron has two state variables: a voltage variable which displays rapid spiking and relaxation under an external stimulus and a slower recovery variable. Each FHN node is governed by the pair of equations
\begin{equation}
\begin{split}
\label{eqn: FHN}
\dot{v}_j & = \alpha_3v_j^3+\alpha_2v_j^2+ \alpha_1v_j - w_j + \frac{K}{n}\sum_{i=1}^n A_{ij}(v_j - v_i) \\
\dot{w}_j & = cv_j - bw_j + \frac{K}{n}\sum_{i=1}^n A_{ij}(w_j - w_i).
\end{split}
\end{equation}
For these simulations we use parameters $\alpha = (-0.1,1.1,-1)$, $c = 0.1$, and $b =0.1$ as studied in \cite{shlizerman2012neural}. This combination yields semicoherent spiking behavior across the network (Fig. \ref{fig: coarsegraining}d). In the raster plot of the network, we see that subsets of the nodes are firing together. The limit-cycle that emerges from projecting this onto 3 SVD modes has a distinctive triangular shape. Each edge appears to result from the spiking of a different subset of nodes. Note that for different parameter values, the system may cease spiking and settle to a steady state or spiking may be more chaotic. For this system we projected to low dimensions by taking the SVD of the voltage and the recovery variables within the nodes separately. $U_1, U_2$ and $U_3$ are the first three modes in the voltage variable. $U_4$, $U_5$, and $U_6$ are the first three modes in the recovery variable. Applying SINDy with trimming to this highly nonlinear, 6-dimensional limit cycle yields a model that is linear in the recovery variables and cubic in the voltage variables. The model predictions and the true dynamics of the system match up very well even over long time horizons  (Fig. \ref{fig: NetworksSINDy}d).  

\begin{figure}[t]
    \centering
    \includegraphics[width=\textwidth]{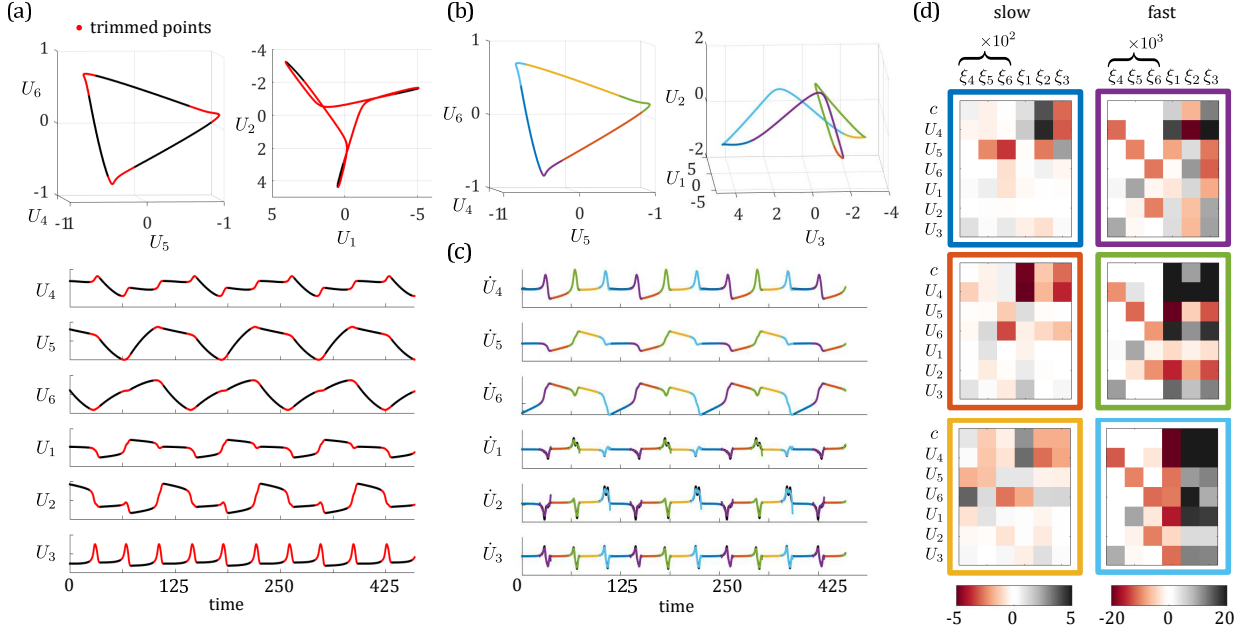}
    \caption{SINDy with trimming can be used to build a linear hybrid model for the low-rank dynamics of an FHN network. (a) Trimming identifies regions of rapid change along directions $U_1$ through $U_4$, which corresponds to the corners of the triangular limit-cycle in the recovery modes $U_6$ and $U_5$. (b) The alternating segments of fast and slow scale dynamics along the limit cycle are used to partition phase-space into 6 regions with unique dominant-balance dynamics. (c) The derivative of each low-rank state variable in each of the six regions are fit using a linear library. The slow dynamic regions (blue, orange and yellow) show great agreement in all variables. The fast dynamic regions are also well matched in modes the recovery modes, 4,5 and 6. The fit is not as clean in the voltage modes, 1, 2 and 3. (d) The coefficients for the 6 dynamic models can be used to assess the system behavior within each region.}
    \label{fig: FHN}
\end{figure}

The multiscale nature of the limit-cycle in the low-rank representation of our FHN network makes it a prime candidate for building a hybrid model. After the network dynamics are reduced to 6 dimensions, we can apply the same steps as outlined for the Rayleigh and Van der Pol oscillators. First, we apply SINDy with trimming. Trimming identifies regions of fast-scale dynamics along directions $U_1$ through $U_4$, which correspond to the three subsets of partially synchronized oscillators firing. In the low-rank dynamics these are the three corners of the triangular limit-cycle on the recovery modes and the three edges of the limit-cycle on the voltage modes (Fig. \ref{fig: FHN}a). Along the limit cycle, alternating segments of fast and slow scale dynamics are used to partition phase-space into 6 regions, each of which has unique dominant-balance dynamics. So rather than just $X_{fast}$ and $X_{slow}$, we have $X^1_{fast}$, $X^2_{fast}$, $X^3_{fast}$, $X^1_{slow}$, $X^2_{slow}$, and $X^3_{slow}$. These six intervals are color coded in Fig. (\ref{fig: FHN}b). The derivatives of the low-rank state variables in each of the six regions are fit using a linear library. We can see that the slow dynamic regions (blue, orange and yellow) show great agreement in all variables. The fast dynamic regions are also well matched in the recovery modes, 4,5 and 6. The fit is not as clean at the ends of each fast region in the voltage modes, 1, 2 and 3. The coefficients for the 6 dynamic models can be used to assess the system behavior within each region (\ref{fig: FHN}c).

\subsection{Heterogeneous Oscillator Networks}
\begin{figure}[htp]
\includegraphics[width=\textwidth]{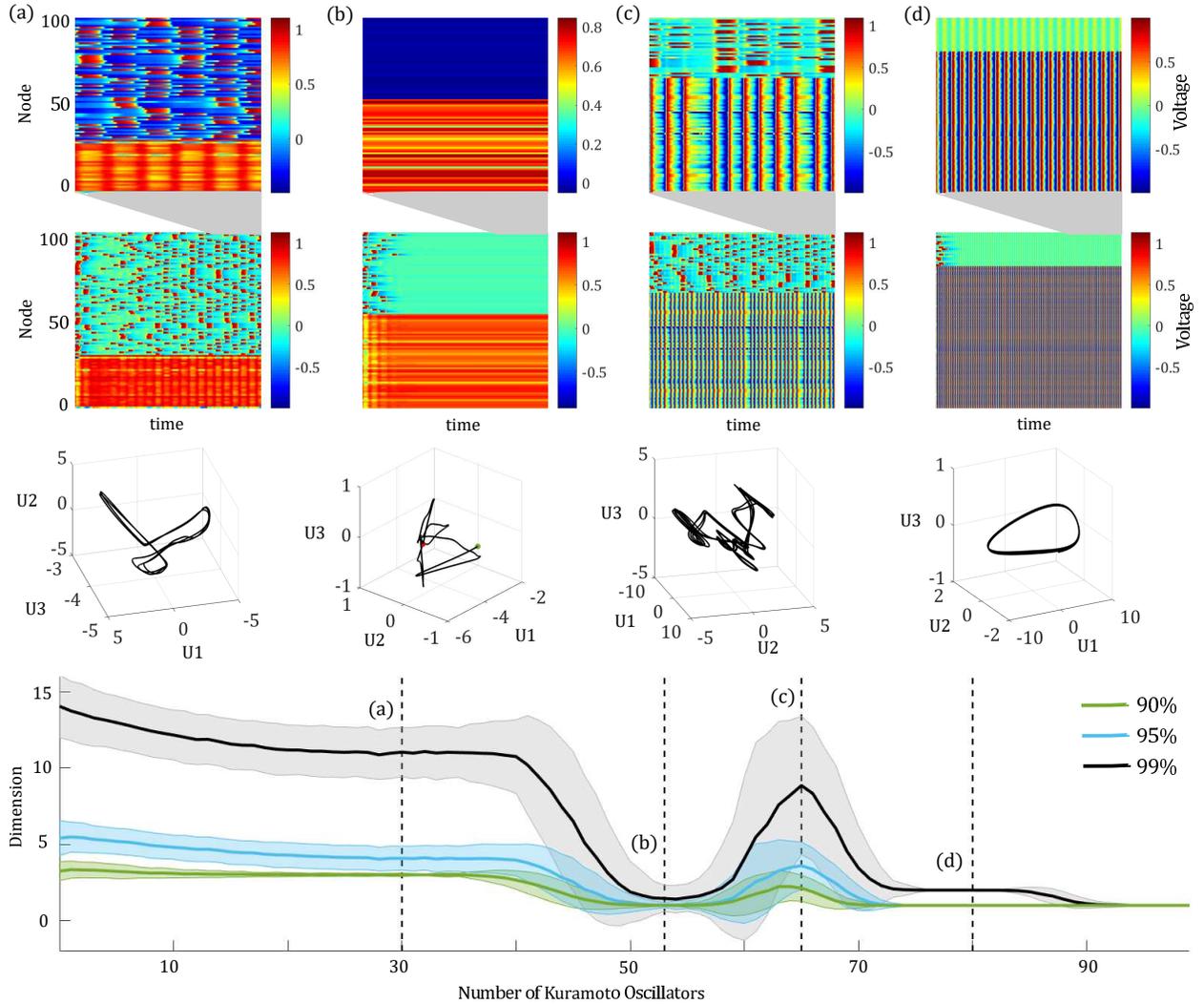}
\caption{Systems of coupled heterogenous oscillators also display low-rank dynamics which can be coarse-grained via the SVD. Here we illustrate dynamics of networks mixing Kuramoto and FHN oscillators in different ratios. Networks without Kuramoto oscillators resemble the FHN behavior shown in Fig. \ref{fig: coarsegraining}, whereas networks with 100\% Kuramoto oscillators have the same exact behavior as the Kuramoto system shown in Fig. \ref{fig: coarsegraining}. When the number of Kuramoto oscillators increases above zero, the Kuramoto oscillators synchronize to what is typically the 3rd FHN mode. This causes the system behavior to develop a progressively stronger pulse along one dimension compared to the purely FHN-like cycle. Eventually, around a 50-50 Kuramoto-FHN split, the majority of simulations progress towards a stable steady state rather than sustaining oscillations. Then as the number of Kuramoto oscillators increases again, there is an interesting region where sometimes the system dynamics resemble those shown in (c) with a dominant oscillatory first mode and the echo of the FHN triangular cycle in modes 2 and 3. Other times the behavior lands on a figure-8 like attractor. As the number of Kuramoto oscillators continues to increase past 70, their combined signal dominates, inducing synchronous low amplitude oscillations in the FHN nodes. One thousand simulations were run for each ratio of Kuramoto to FHN oscillators and the dimension of the resulting dynamics was estimated by tracking the number of modes needed to capture 90, 95, and $99\%$ of the Frobenius norm. The mean ``dimensions" as defined by each threshold are shown by solid lines and the corresponding standard deviations are shaded.}
\label{fig: 2types}
\end{figure}

\begin{figure}[t]
\includegraphics[width=0.9\textwidth]{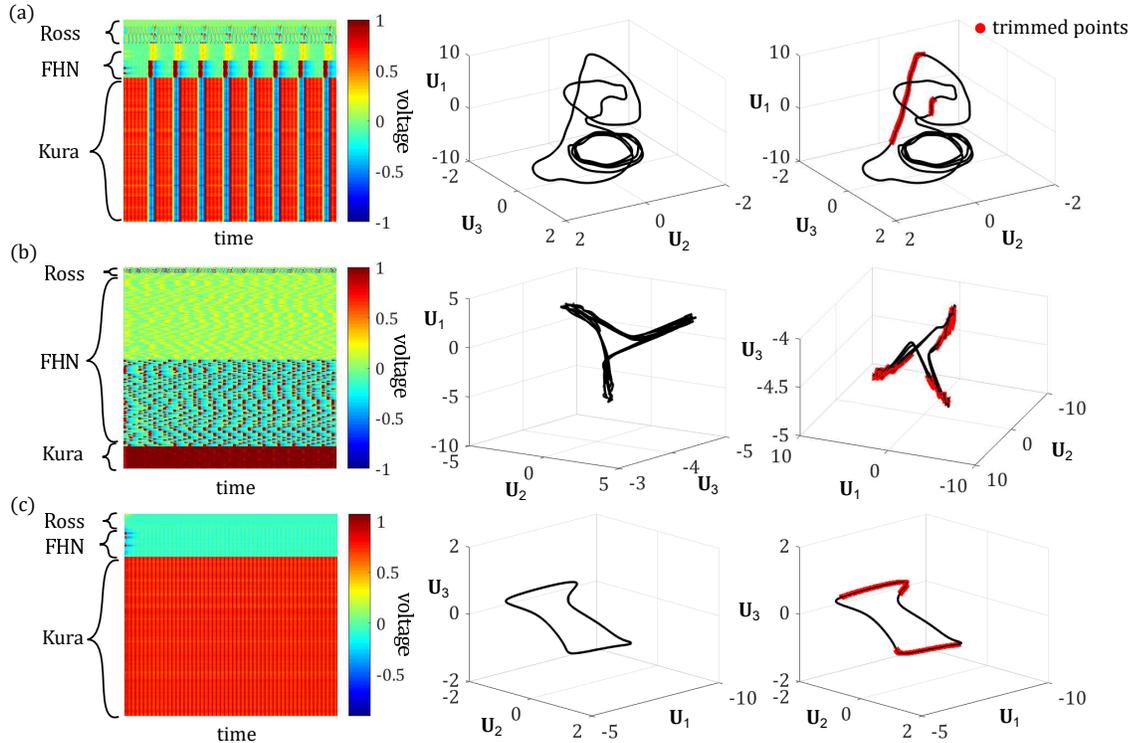}
\caption{Heterogenous oscillator networks mixing Kuramoto, FHN, and Rossler nodes can also display low rank dynamics in certain regimes. Here we show several examples and the corresponding limit cycle resulting from projecting the system onto its first 3 singular vectors. SINDy with trimming again partitions the trajectory into fast vs. slow-scale dynamics. (a) All nodes across the system synchronize to a background frequency with low amplitude, on top of which there are periodic large spikes. Trimming pulls out the rapid relaxation following these large spikes and the steepest portion of the spike. (b) When the majority of the nodes are FHN, the system is dominated by FHN-like behavior. However, for this system there are small spikes along the three corners of the FHN-cycle resulting from the Rossler nodes. Trimming identifies the contribution of the Rossler nodes to the overall limit cycle. (c) Here Kuramoto oscillators dominate system behavior, but oscillations are damped because the FHN nodes have settled to a steady state. This produces an elegant non-linear limit cycle. Trimming pulls out regions of rapid change along $U_1$.}
\label{fig: 3types}
\end{figure}

By considering networks of nodes displaying relaxation-oscillations and chaotic dynamics, we demonstrated the power of this method to disambiguate and capture spatio-temporal heterogeneity. However, we are not limited to networks consisting of a single type of oscillator. We can enhance the heterogeneity of the system by mixing two or three types of oscillators together in one network and apply the same coarse-graining method (Fig. \ref{fig: 2types} and \ref{fig: 3types}). 

As a case study on networks of two types of oscillators, we coupled together $n_K$ nodes governed by Kuramoto dynamics and $n_F$ nodes governed by FHN dynamics where $n_K + n_F = 100$. When coupling FHN nodes to Kuramoto nodes equations (\ref{eqn: Kuramoto}) and (\ref{eqn: FHN}) apply to the appropriate subset of nodes in the system, except that FHN nodes are coupled to the sine of the Kuramoto node to maintain a consistent scale. For visualization purposes we ordered the system such that the first $1,.., n_k$ nodes are Kuramoto and the next $n_k+1, ..., 100$ are FHN. The adjacency matrices for these simulations are Erdos-Renyi graphs with connection probability $p = 0.2$. The FHN parameters are $\alpha = (-0.1,1.1,-1)$, $c = 0.1$, and $b =0.1$ with connectivity strength $K_F = 0.2$. The intrinsic frequencies for Kuramoto oscillators are drawn from a uniform distribution centered at $0.6$, so the population would naturally oscillate much more quickly than the FHN firing rate. The connectivity strength for Kuramoto oscillators is $K_K = 10$, so this bloc also synchronizes more quickly than the FHN bloc. 

With this set of parameters, we can simulate many different network instances at all possible ratios of Kuramoto to FHN oscillators. When $n_K = 0$, the system matches the all FHN network discussed previously with collective dynamics landing on a triangular limit cycle. When $n_K = 100$, the system matches the all Kuramoto network, with low rank dynamics attracted to a stable, elliptical limit cycle. For ratios in between the low-rank dynamics evolve from the triangular FHN cycle to the elliptical Kuramoto cycle (Fig. \ref{fig: 2types}). When $n_K \in [1,40]$, the FHN bloc dominates the collective dynamics and the Kuramoto oscillators synchronize to what is typically the third mode of the FHN network. As the number of Kuramoto oscillators increases, this strengthens the influence of the third mode, widening the variation in this direction and rounding out the sharp corners of the triangular cycle. When the numbers of Kuramoto and FHN nodes are more equal, around $n_K \in [40, 60]$, the system dynamics often collapse to a stationary state over time. Then continuing to increase $n_K$ into the range $[60,70]$, there is a window of interesting mixed dynamics. In this range, the Kuramoto bloc sets the pace for a background pulse that is shared across all nodes, including the FHN nodes. However, the FHN nodes still fire periodically. The corresponding low rank dynamics always have a dominant oscillatory first mode. For some networks, the FHN recovery periods stretch the period of the Kuramoto oscillators, and so an echo of the FHN triangle is captured in modes 2 and 3 as shown in Fig. \ref{fig: 2types}c. In other cases, the Kuramoto oscillators are not strongly entrained to the FHN firing, so the second mode is also oscillatory and mode 3 captures chaotic oscillations. These collective dynamics land on a figure-8 shaped attractor. Finally, increasing $n_K$ still further into the range $[70, 90]$ causes the Kuramoto bloc to dominates as it entrains the FHN nodes to pulse at the mean Kuramoto frequency. These synchronized oscillations are well captured by a curved ellipse in the SVD projection. The collective low-rank dynamics appear almost entirely Kuramoto-like, ie. follow an elliptical limit cycle, beyond $n_k = 90$. In this example, we thoroughly explored the dynamics observed for one set of parameters on a mixed Kuramoto-FHN network. However, many other potential parameter sets could be considered and would lead to different collective dynamics. 

\begin{figure}[htp]
    \centering
    \includegraphics[width=\textwidth]{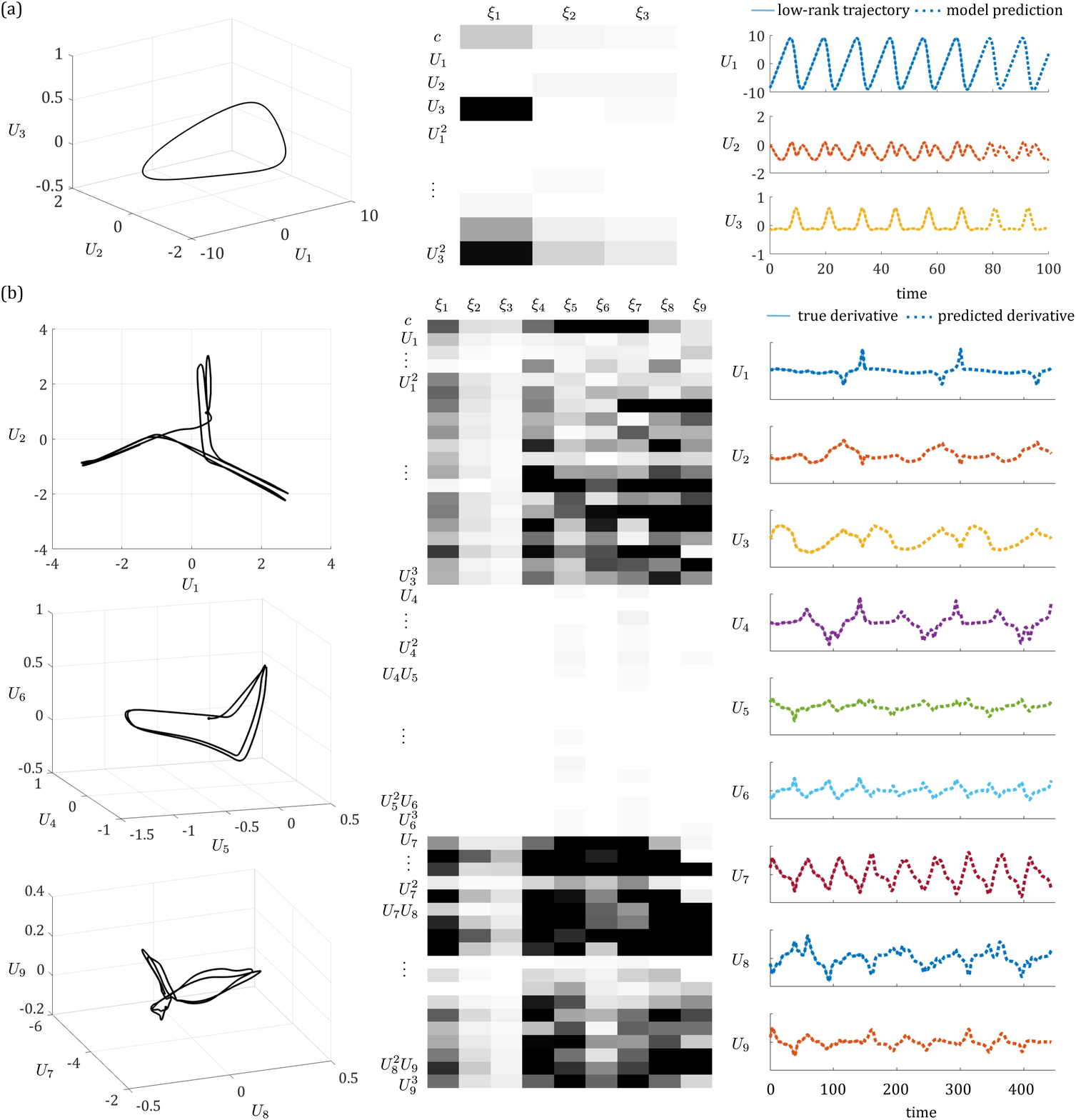}
    \caption{SINDy can be applied to the low-rank representation of collective dynamics of mixed FHN-Kuramoto networks. (a) In cases where the collective dynamics are quite low-dimensional and relatively simple, a stable model is recovered. (b) However, for a network with higher intrinsic dimension to its dynamics that is navigating a highly non-linear limit cycle, we may need to use many modes to accurately match the derivative and the resulting models may be quite dense. }
    \label{fig: 2typesSINDy}
\end{figure}

For each ratio of $n_K$ to $n_F$ and this fixed set of parameters, we estimated the dimensions of the collective dynamics by calculating the number of SVD modes required to reconstruct the network dynamics to $90$, $95$, and $99\%$ accuracy with respect to the Frobenius norm. This estimate was repeated for 1000 network iterations at each ratio, and the mean and standard deviation were calculated (Fig. \ref{fig: 2types}). The expected dynamic dimension of the network is the highest when the FHN bloc is dominant. It is low around a 50:50 ratio when the system often  goes to a steady state and again above $n_K = 90$ when the Kuramoto bloc fully dominates. There is a jump in dimension for the window around $n_K = 60$ to $70$ where we see the system exhibiting collective dynamics that are a fairly even balance between the two styles of oscillator. These estimates of the dynamic dimension give a hint as to how many modes need to be used when fitting a SINDy model to the collective limit cycle. In Fig. \ref{fig: 2typesSINDy}, we apply our coarse graining method to produce models for the example in panels (a) and (d) of Fig. \ref{fig: 2types}. For the network in panel Fig. \ref{fig: 2typesSINDy}a, which is mostly dominated by Kuramoto oscillators, though the FHN influence makes the limit cycle non-linear, our method produces a stable model that can be used to predict network dynamics over a long time period. The first mode is quite dominant, we see much larger coefficients on these dynamics compared to the other two modes. In contrast, to accurately fit the derivative for the network that is FHN dominated, we need 9 modes and a cubic library. In this case the SVD was applied to the Kuramoto and the FHN blocs separately, akin to learning coarse grained variables for subsets of the nodes in a network as in the method of Antonsen and Ott  \cite{ott2009long}. The voltage variables of the FHN nodes and the Kuramoto operate on a fast scale, whereas the recovery variables of the FHN nodes operate on a slow scale. The recovered model fits the derivative well, however it is unstable. This is an inherent challenge in the standard SINDy implementation, which can be addressed through further modifications that ensure model stability. There are several stabilization approaches introduced in \cite{delahunt2021toolkit, hirsh2021sparsifying}. 

\begin{figure}[htp]
\includegraphics[width=\textwidth]{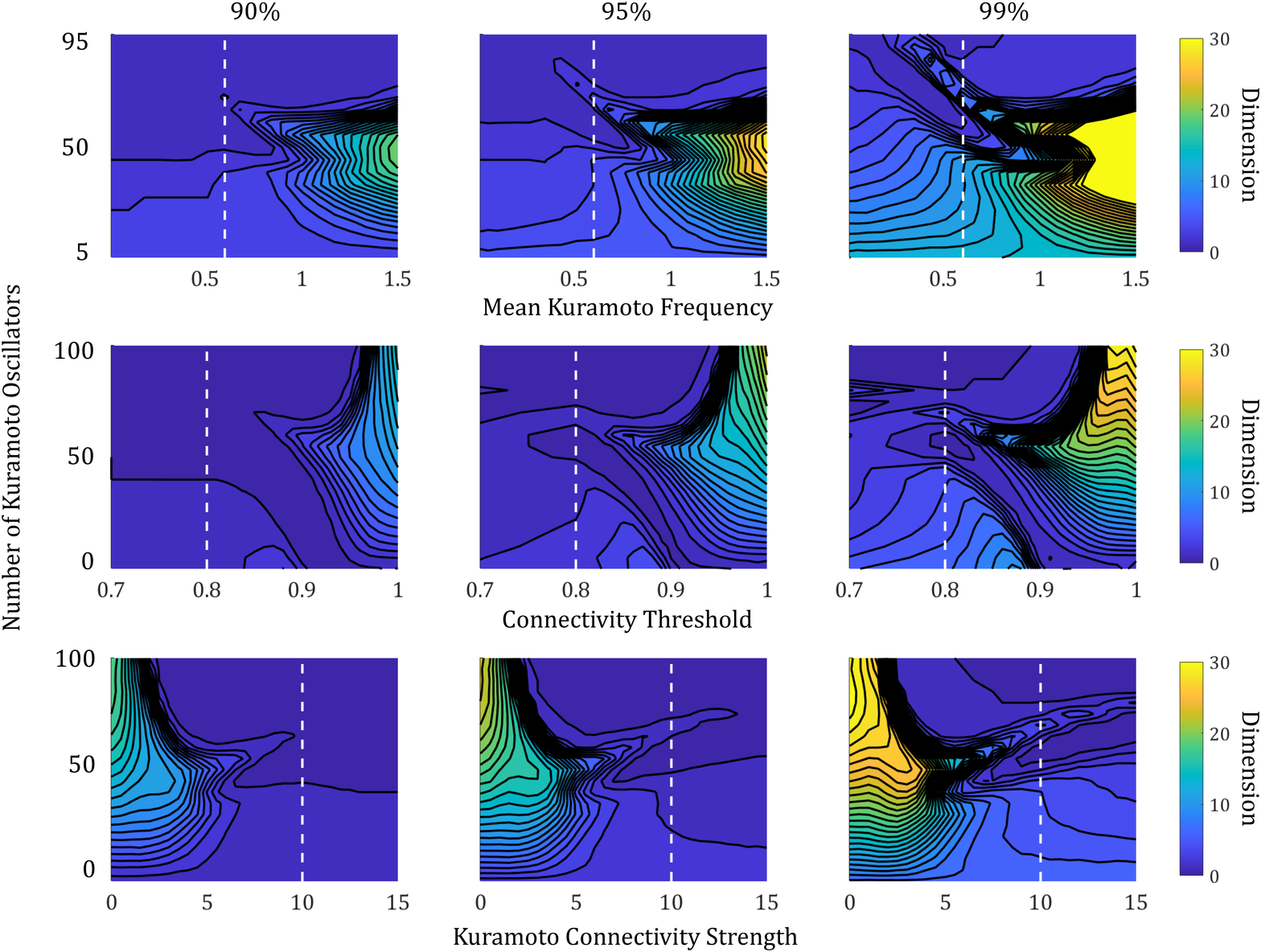}
\caption{The typical dimension of dynamics for systems of coupled Kuramoto and FHN oscillators is a function of several system parameters: the ratio of Kuramoto to FHN oscillators in the network, the mean intrinsic frequency of Kuramoto oscillators, the Erdos-Renyi connectivity threshold used to construct the adjacency graph, and the connectivity strength among the Kuramoto oscillators. For each parameter combination, one thousand independent systems were simulated and the mean number of modes needed to capture $90\%$, $95\%$, and $99\%$ of the Frobenius norm was calculated. In row 1, observe that as the mean frequency of the Kuramoto oscillators increases, the dimension of systems with substantial numbers of both Kuramoto and FHN nodes increases. However, systems that are dominated by Kuramoto-type behavior remain low dimensional across all frequencies. In row 2, note that as the connectivity threshold increases (coupling is universally weaker) the Kuramoto dominant networks rapidly jump in dimension reflecting the bifurcation between synchronicity and asynchronicity in the system. On the other hand, the FHN-dominant networks actually decrease in dimension as coupling weakens because without sufficient input from neighbors, these nodes will cease spiking. The third row shows variation in the connectivity specifically within the Kuramoto population. In this case higher connectivity strength means a more connected Kuramoto bloc. Again, as connectivity weakens the dimension of Kuramoto-dominant systems increases sharply below a threshold around 2.5. For networks with small numbers of Kuramoto oscillators, the system dimension stays fairly constant as the connectivity strength changes because the FHN bloc is dominant and remains unaffected. The white dotted line indicates the parameter regime explored in Fig. \ref{fig: 2types}.}

\label{fig: ParameterMaps}
\end{figure}

In addition to varying the ratio of Kuramoto to FHN oscillators, we also estimated the dimension of collective network dynamics as a function of other network parameters.  We tested a wide range of network connectivities and differed the scale between the FHN and Kuramoto oscillations. These explorations are mapped out in Figure \ref{fig: ParameterMaps}. Note that the white dotted line indicates the parameter set that was illustrated in Fig. \ref{fig: 2types}: the mean Kuramoto frequency is 0.6, the connectivity threshold (1 minus the connectivity probability) on the adjacency matrix is $0.8$, and the Kuramoto connectivity strength is 10, while $n_K \in [0, 100]$. When constructing these figures, the dimension of a given instance of a network was estimated by simulating the network until time $t = 2000$, taking the SVD of the system dynamics from time $t = 1000$ onwards, and then calculating the number of singular values required to retain $90, 95$ and $99\%$ accuracy in reconstruction. Imposing a time delay before computing the SVD allows the system to settle onto the long-term limit cycle, so our estimate of the dimension discounts the early transient behavior of the system. For a given parameter set, 1000 instances of independent networks were simulated and the average dimension for each accuracy threshold was calculated. As we change two parameters defining the network structure and two parameters defining the Kuramoto nodes' internal dynamics, we observe wide variation in the typical dimension of the collective dynamics. 

Starting in the top row of Figure \ref{fig: ParameterMaps}, we see that the dimension of the collective dynamics increases sharply with the mean Kuramoto frequency when the ratio $n_K:n_F$ is around 1. At higher frequencies, the intrinsic Kuramoto dynamics are on a very different time scale than the instrinsic FHN dynamics so we no longer achieve entrainment of one group by the other for equal size blocs. There is a narrow region of low dimension that persists for oscillator ratios above 50 up to a mean frequency of about 1. This is the region in which many system instances go to a steady stationary state. When the Kuramoto bloc is dominant, the dimension does not vary significantly as the mean frequency changes because the oscillators synchronize to one speed; no matter how fast or slow that speed is, the dynamics will be low-dimensional. When the FHN bloc is fully dominant, the dimension does not vary significantly with the kuramoto oscillator frequency because the dynamics of all nodes are driven almost exclusively by coupling to the FHN nodes.

In the second row of Figure \ref{fig: ParameterMaps}, we vary the connection probability, $p$, when constructing the adjacency matrix for the network. As the connectivity threshold, $1-p$, increases, the coupling across all nodes is weakened until the network becomes disconnected. In the upper half of these maps, when the collective dynamics are dominated by Kuramoto oscillators, we observe the phase transition from a synchronized or semi-sychronized state with dimension less than 3 to an asynchronous state with dimensions well above 10. This transition is apparent as a black region for high $n_K$ and high connectivity thresholds because the dimension is increasing so rapidly that the contour lines are on top of each other. When FHN nodes make up a significant amount of the population, this rapid jump in dimension is ameliorated because as FHN nodes become disconnected, they lose the driving force needed to maintain tonic spiking instead relaxing to a stable steady state. The constant signal from the FHN nodes reduces the overall dimension of the system. When FHN nodes dominate the system ($n_K < 20$), there is a phase transition from tonic spiking behavior for connectivity thresholds in the range $[0.7, 0.9)$ to a constant steady state above $0.9$, again visible as a black region where contour lines are on top of each other in the $99\%$ plot.

Finally, we examined the dimension of network dynamics as we varied the connectivity strength amongst the bloc of Kuramoto oscillators (Fig. \ref{fig: ParameterMaps}, row 3). In this case higher values of $K_K$ indicate stronger coupling and higher synchronization within the Kuramoto bloc. Lower connectivity leads to higher dimension collective dynamics and higher connectivity leads to lower dimension collective dynamics if the number of Kuramoto oscillators is above 1. As with the connectivity threshold, we observe a sharp phase transition between partial synchronicity and asynchronicity as a black region on the maps when $n_K > 50$. For smaller values of $n_K$ the transition is more gradual because the FHN modes are well-connected and offset the chaotic firing of the Kuramoto bloc. When $n_K$ is very small, the dynamic dimension is independent of $K_K$ because the FHN dynamics are fully dominant.

We also applied our coarse-graining method to networks consisting of three different types of oscillators: Kuramoto, FHN, and Rossler. With this level of heterogeneity, the possible sets of network parameters to consider quickly grows, making a thorough exploration intractable. However, for some parameter regimes and initial conditions we demonstrate that the collective dynamics are low-dimensional and can be captured by a small number of SVD modes. Possible behaviors include high frequency, low amplitude background oscillations with periodic high amplitude spikes across the full system (Fig. \ref{fig: 3types}a), FHN-like semi-synchronized spiking with Rossler oscillators adding higher-frequency oscillations during the relaxation period (Fig. \ref{fig: 3types}b), or a repeated high-frequency, nonlinear oscillation at a low amplitude across the Kuramoto bloc (Fig. \ref{fig: 3types}c). In each of these examples, we can apply SINDy with trimming to isolate regions of differing dynamics. In panel (a) from Fig. \ref{fig: 3types}, the points corresponding to the universal rapid relaxation across the system are trimmed along with the steepest portion of the collective spike. In panel (b) from Fig. \ref{fig: 3types}, the low dimensional dynamics look like the now familiar triangular FHN limit cycle, but there are irregular oscillations allow the corner wings of the triangle. There portions seem to be contributed by the small bloc of Rossler oscillators in the system and are precisely the points picked out by trimming. In panel (c) from Fig. \ref{fig: 3types}, trimming identifies regions of rapid change along the $U_1$ dimension, while $U_3$ is constant and $U_2$ is nearly constant.

\section{Discussion}

Modeling heterogeneous, multiscale physics remains a mathematically challenging proposition.  Indeed, traditional computational techniques quickly become intractable when attempting to resolve such systems in both space and time.  Emerging data-driven methods are enabling new mathematical architectures that can automate the discovery of coarse-grained coordinates and dynamics that are low-dimensional and parsimonious.  This renders approximate models that allow for rapid evaluation of the heterogeneous, multiscale physics.  More than that, these techniques can also be constructed to produce interpretable physics models.

We have proposed a set of algorithms for modeling heterogeneous, high-dimensional networked dynamical systems.  Specifically, we have leveraged dimensionality-reduction, sparse regression, and robust statistics to discover interpretable, coarse-grained models of the underlying physics.  The method further allows us to disambiguate between the resulting fast and slow scale dynamics of the system. Used in combination, these methods are able to: (i) identify low-dimensional embeddings (coordinates) on which to construct models, (ii) identify different time scales using the statistical robustification technique of data trimming, and (iii) identify interpretable parsimonious models of the coarse-grained dynamics at different timescales.

The mathematical architecture has been demonstrated with a series of numerical experiments on networked nonlinear oscillators.  The oscillators, whether a single type of oscillator or a heterogeneous set of oscillators, evolve and interact to produce low-rank dynamics that often include relaxation oscillations. The literature is largely devoid of the consideration of such systems due to their complexity.  However, the mathematical architecture developed here is able to extract meaningful models even with such high-dimensional heterogeneous interactions.

\section*{Acknowledgements}

The authors acknowledge funding support from the Air Force Office of Scientific Research (AFOSR FA9550-19-1-0386) and from the National Science Foundation AI Institute in Dynamic Systems grant number 2112085.

\medskip

\bibliographystyle{unsrt}
\bibliography{main}

\end{document}